\documentclass[12pt]{article}
\usepackage[utf8]{inputenc}
\usepackage{amsmath, amssymb, amsfonts}
\usepackage{amsthm}
\usepackage{graphicx}
\usepackage[dvipsnames,svgnames,x11names]{xcolor}
\usepackage{fancyhdr}
\usepackage[sort&compress,comma,square,numbers]{natbib}
\usepackage[top=2cm,bottom=2cm,left=2cm,right=2cm]{geometry}
\usepackage{multicol}
\usepackage{booktabs}

\theoremstyle{plain}
\newtheorem{thm}{Theorem}[section]

\newtheorem{corollary}[thm]{Corollary}
\newtheorem{lemma}[thm]{Lemma}

\theoremstyle{definition}

\newtheorem{remark}[thm]{Remark}

\begin{document}

\section*{INTRODUCING AND APPLYING S.C.E MODEL UNDER DUSART'S INEQUALITY TO PROVE GOLDBACH'S STRONG CONJECTURE FOR 74 TYPICAL STRUCTURES OUT OF ALL 75 STRUCTURAL TYPES OF EVEN NUMBERS}

\begin{center}
AREF ZADEHGOL MOHAMMADI${ }^{1,2}$, MOHSEN KOLAHDOUZ${ }^{3,4}$
\end{center}

\subsection*{Abstract}

In this paper, we present a relative proof for Goldbach's strong conjecture. To this end, we first present a heuristic model for representing even numbers called Semi-continuous Model for Even Numbers or briefly S.C.E Model, and then by using this model we categorize all even numbers into 75 distinct typical structures. Also in this direction, we employ this model along with the following inequality to obtain the relative proof

\begin{equation}
\frac{x}{\ln x} \leq_{x \geq 17} \pi(x) \leq_{x>1} 1.2251 \frac{x}{\ln x}
\end{equation}

where $\pi(x)$ denotes the number of all primes smaller than and equal to $x$. This inequality is presented by Pierre Dusart in his paper [P. Dusart, Explicit estimates of some functions over primes, Ramanujan J. 45 (2016), No. 1, 227-251].

In fact, by relative proof we mean that 74 typical structures out of 75 ones satisfy Goldbach's strong conjecture. Also, since the last typical structure is the dominant structure over all even numbers, we come up with three unproven inequalities for elements of S.C.E model using each of which, we can prove Goldbach's strong conjecture for this structure too. It is necessary to say that, we guess theses three inequalities can be proved the same as to Dusart's inequality.     

\noindent 2010 Mathematics Subject Classification: Primary 11P32, 11A41; Secondary 11A67, 11N05.

\noindent \textbf{Keywords:} Goldbach's Strong conjecture, Semi-Continuous Model of even numbers, Prime numbers, Even numbers, Odd numbers, Dusart's inequality, Additive interaction of even numbers

\section{History}

Basically, the most prominent narration to the content of Goldbach's strong conjecture refers to the possibility of converting any even number greater than or equal to 4 to sum of two prime numbers; which was expressed by Christian Goldbach in 1742. For years, there have been told many alternatives to this conjecture; but in overall, this conjecture is separated into two branches of weak and strong types, which the latter case refers to the possibility of converting any odd number greater than 5 into the sum of three prime numbers. The proof of weak type conjecture is widely accepted via a publication presented by Harald Helfgott Anderson in 2013[1].

However, the strong type conjecture is still under examination and investigation of mathematicians in order to find an argument to prove the generality of it or to find a counterexample to revoke its generality.

\section{Introduction and Preliminaries}

In the process of studying and working on even numerical examples in order to find out why or how the process of Goldbach's conjecture holds true, it can be seen that all even numbers can be considered as a distance from zero with even values. Due to this point of view, in this paper, we will observe a heuristic model in correspondence with the Goldbach's strong conjecture that we call it Semi-Continuum Model of Even Numbers or briefly S.C.E Model.

By using this model, we can translate Goldbach's strong conjecture into this model by focusing on $d_{E}$ element which will be introduced later. It is necessary to mention that we prove this conjecture relatively just for even numbers $E>2525$. On the other hand, by using computer's application or manually, it can be easily seen that for all even numbers $E<2525$ the Goldbach's strong conjecture is valid. Therefore, we can extend our findings to all even numbers.

Since in this model only additive interactions of odd numbers are involved, the number 2 despite being a prime will not be considered. However, by expressing Goldbach's strong conjecture via this model, we can equivalently say that every even number $E \geq 4$ can be viewed as a sum of two primes if and only if the element $d_{E}$ never vanishes.

In this direction, we employ also the precious inequality presented by Pierre Dusart in [2] as

\begin{equation}
\frac{x}{\ln x} \leq_{x \geq 17} \pi(x) \leq_{x>1} 1.2251 \frac{x}{\ln x}
\end{equation}

where $\pi(x)$ denotes the number of all primes smaller than and equal to $x$. Here, the first inequality holds for $x \geq 17$, and the second one holds for $x>1$. The interested reader can observe other inequalities for $\pi(x)$ in that paper, but here the inequality (1) is considered for both simplicity and universality. It will be turned out that by linking or importing this inequality into S.C.E Model, we can come up with the arguments, of which yields that for all even numbers $E \geq 4$, the element $d_{E}$ never vanishes, at least for 72 typical structures.

\subsection{Basic Semi-continuum Model of Even Numbers or S.C.E Model}

In this section, for each positive even number $E$ we uniquely correspond a quadruple representation as

\begin{equation}
\frac{E}{4}=a_{E}+b_{E}+c_{E}+d_{E}
\end{equation}

where $a_{E}, d_{E} \in\left\{\frac{n}{2} \mid n \in \mathbb{N}\right\}$ and $b_{E}, c_{E} \in \mathbb{N}$. To describe these elements, we first define the notion of additive interaction of odd numbers for a given even number $E$ as follows.

\textbf{Definition 2.1.} Let $E$ be a given even number, then for all odd positive numbers $x, y<E$ the notion

\begin{equation}
x \sim y
\end{equation}

is called an additive interaction for $E$ provided $x \leq y$ and $x+y=E$.

In this setting, for all additive interactions of $E$, we determine the elements of (2) as

\begin{align}
a_{E} &= \sharp\{x \sim y \mid x \neq y \text{ and } x, y \text{ are nonprimes}\} + \frac{1}{2} \sharp\{x \sim x \mid x \text{ is nonprime}\} \\
b_{E} &= \sharp\{x \sim y \mid x \text{ is nonprime but } y \text{ is prime}\} \\
c_{E} &= \sharp\{x \sim y \mid x \text{ is prime but } y \text{ is nonprime}\} \\
d_{E} &= \sharp\{x \sim y \mid x \neq y \text{ and } x, y \text{ are primes}\} + \frac{1}{2} \sharp\{x \sim x \mid x \text{ is prime}\}
\end{align}

where $\sharp$ denotes cardinality symbol.

\textbf{Definition 2.2 (Semi-continuum model of even numbers).} Let $E$ be an even number, then the unique quadruple representaion (2) is called the semi-continuum model of $E$, where $a_{E}, b_{E}, c_{E}, d_{E}$ are as (3)-(6), respectively.

In this case by defining $L_{1}(E)=a_{E}+b_{E}$, $L_{2}(E)=c_{E}+d_{E}$, $R_{1}(E)=a_{E}+c_{E}$ and $R_{2}(E)=b_{E}+d_{E}$, we can also demonstrate the semi-continuum model of even number $E$ as follows.

\begin{center}
\includegraphics[width=0.5\textwidth]{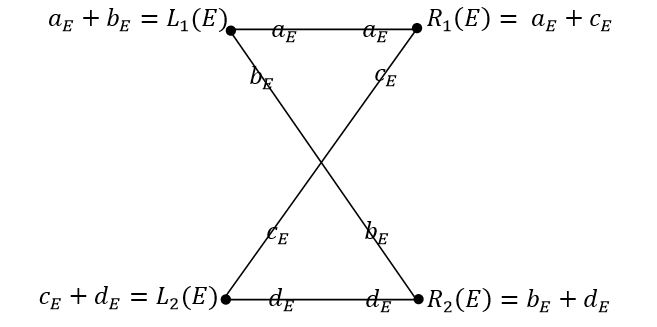}
\end{center}

\begin{center}
Figure 1. S.C.E Model of $E$
\end{center}

In this papion shape, for upper vertices, we can see that the values

\begin{equation}
\left\lceil L_{1}(E)\right\rceil,\left\lceil R_{1}(E)\right\rceil
\end{equation}

show respectively the number of odd non-primes in the intervals $\left[0, \frac{E}{2}\right]$ and $\left[\frac{E}{2}, E\right]$, where $\lceil \cdot \rceil$ is the ceiling function. Similarly, for lower vertices, the same reasoning yields the values

\begin{equation}
\left\lceil L_{2}(E)\right\rceil,\left\lceil R_{2}(E)\right\rceil
\end{equation}

show the number of primes in the intervals $\left[0, \frac{E}{2}\right]$ and $\left[\frac{E}{2}, E\right]$, respectively. Furthermore, by (2), we can also drive the following relation

\begin{equation}
\frac{E}{4}=L_{1}(E)+L_{2}(E)=R_{1}(E)+R_{2}(E)
\end{equation}

As it is obvious, since only the odd numbers $x \leq y<E$ additively interact with each other to represent the even number $E$, this model is named for semi-continuum model of even numbers.

\textbf{Remark 2.3.} We can easily see that

\begin{equation}
b_{E}+c_{E}+2 d_{E}=L_{2}(E)+R_{2}(E)=\pi(E)-1
\end{equation}

Note that, since 2 is a prime but is not an odd number, we subtract 1 from $\pi(E)$ to get $L_{2}(E)+R_{2}(E)$. On the other hand, we have

\begin{equation}
b_{E}+c_{E}+2 a_{E}=L_{1}(E)+R_{1}(E)=\frac{E}{2}-\pi(E)+1
\end{equation}

\textbf{Remark 2.4.} In view of Remark 2.3, the even number 4 is the only exception to our model since $d_{4}=0$ but we can write $2+2=4$.

\textbf{Remark 2.5.} For an even number $E$, if $\frac{E}{2}$ is also an even number, then we have

\begin{align}
& L_{1}\left(\frac{E}{2}\right)+R_{1}\left(\frac{E}{2}\right)=L_{1}(E) \\
& L_{2}\left(\frac{E}{2}\right)+R_{2}\left(\frac{E}{2}\right)=L_{2}(E)
\end{align}

and if $\frac{E}{2}$ is an odd number, then we have

\begin{align}
L_{1}\left(\frac{E}{2}-1\right)+R_{1}\left(\frac{E}{2}-1\right)+\frac{1}{2} & =L_{1}(E) \\
L_{2}\left(\frac{E}{2}-1\right)+R_{2}\left(\frac{E}{2}-1\right) & =L_{2}(E)
\end{align}

\textbf{Example 2.6.} Let $E=20$, then by considering the distance of 20 from zero as
\begin{center}
\includegraphics[width=0.5\textwidth]{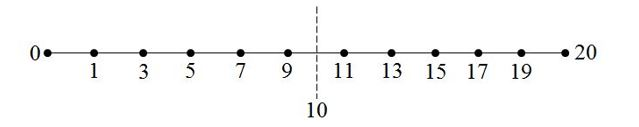}
\end{center}

we can drive all additive interaction of 20 as follows

\begin{center}
\includegraphics[width=0.2\textwidth]{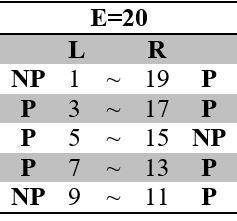}
\end{center}

Now, in view of the representation (2), we can compute quadruple elements of 20 as

\[
a_{20}=0, \quad b_{20}=2, \quad c_{20}=1, \quad d_{20}=2
\]

using which we can obtain the S.C.E Model of 20 as follows

\begin{center}
\includegraphics[width=0.5\textwidth]{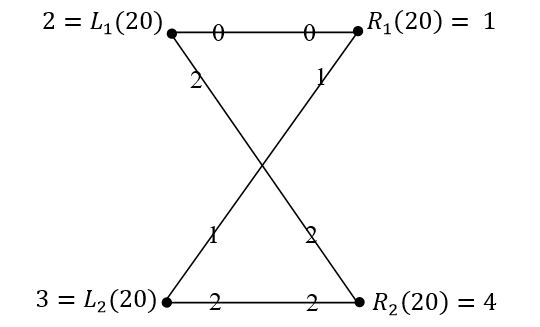}
\end{center}

\begin{center}
Figure 2. S.C.E Model of 20
\end{center}

\subsection*{2.2. Linking with Dusart's Inequality}
Actually, the S.C.E Model by itself is slightly poor to prove Goldbach's strong conjecture. To enrich it, we utilize Dusart's inequality (1). To link up with this inequality, we consider the following two main inequalities called \textit{Teeter Inequalities}.

\textbf{Lemma 2.7.} Let $E \geq 17$ be an even number with quadruple representation
\[
\frac{E}{4}=a_{E}+b_{E}+c_{E}+d_{E}
\]
obtained by (2). Then we have

\begin{align}
\label{33}
& \frac{E}{2}-1.2551 \frac{E}{\ln E}+1<b_{E}+c_{E}+2 a_{E}<\frac{E}{2}-\frac{E}{\ln E}+1 \\
\label{34}
& \frac{E}{\ln E}-1<b_{E}+c_{E}+2 d_{E}<1.2551 \frac{E}{\ln E}-1
\end{align}

\begin{proof}
The assertions of the lemma are respectively direct consequences of relations (8), (9) and Dusart's inequality (1).	
\end{proof} 

\textbf{Corollary 2.8.} By subtracting the inequality (13) from (12), we can conclude the following inequality
\begin{equation}\label{2.1}
	\frac{E}{\ln E}-\frac{E}{4}-1<d_{E}-a_{E}<1.2551 \frac{E}{\ln E}-\frac{E}{4}-1
\end{equation}

\textbf{Remark 2.9.} With the aid of a calculus approach, it turns out that the point $x=$ 130.4574578 is the root of the following decreasing function

\[
f(x)=1.2551 \frac{x}{\ln x}-\frac{x}{4}-1
\]

Thus, in view of (14), for even numbers $E>130.4574578$ we can follow that
\[
d_{E}<a_{E}
\]
	\begin{remark}\label{3.3}
	For all positive even number $E\geq 17$ we have
	\begin{equation}
		\dfrac{E}{4}-1.2551\dfrac{E}{\ln E}+1<a_E.
	\end{equation}
\end{remark}
\begin{proof}
	Obviously, since $ \frac{E}{4}-1.2551\frac{E}{\ln E}+1\notin\{\frac{n}{2}\,\vert\,n\in\mathbb{N}\} $, then $ \frac{E}{4}-1.2551\frac{E}{\ln E}+1\neq a_E $ for all psitive even number $E$. In the sequel, we continue by contradiction. So, let $E$ be a positive even number such that $a_E<\frac{E}{4}-1.2551\frac{E}{\ln E}+1$ or
	\begin{equation}\label{36.1}
		2a_E<\frac{E}{2}-2\frac{1.2551E}{\ln E}+2.
	\end{equation}
	If we subtract the inequality (\ref{36.1}) from (\ref{33}) we get
	$1.2551\frac{E}{\ln E}-1<b_E+c_E$, which contradicts with the inequality (\ref{34}).
\end{proof}	
\begin{corollary}\label{3.4}
	For all positive even number $E\geq22864$ we have
	\begin{equation}\label{37.1}
		b_E+c_E+d_E<a_E.
	\end{equation} 
\end{corollary}
\begin{proof}
	First by the inequality (\ref{34}) we have
	$$b_E+c_E+d_E<1.2551\dfrac{E}{\ln E}-1,$$
	but by calculus we can see 
	$\frac{x}{4}-1.2551\frac{x}{\ln x}+1>1.2551\frac{x}{\ln x}-1$ for all real number $x\geq 22864$, then we can easily obtain (\ref{37.1}) by previous remark.
\end{proof}
To obtain a lower and a upper bound for $L_{1}(E), L_{2}(E), R_{1}(E)$ and $R_{2}(E)$ individually with the aid of Dusart's inequality, we first compute $L_{2}(E)$ and the rest will be readily accessible.

\textbf{Lemma 2.10.} Let $E>34$ be an even number, then the following inequalities hold true
\[
\begin{aligned}
& \frac{\frac{E}{2}-1}{\ln \left(\frac{E}{2}-1\right)}-1<L_{2}(E)<1.2551 \frac{\frac{E}{2}}{\ln \frac{E}{2}}-1 \\
& \frac{E}{4}-1.2551 \frac{\frac{E}{2}}{\ln \frac{E}{2}}+1<L_{1}(E)<\frac{E}{4}-\frac{\frac{E}{2}-1}{\ln \left(\frac{E}{2}-1\right)}+1 \\
& \frac{E}{\ln E}-1.2551 \frac{\frac{E}{2}}{\ln \frac{E}{2}}<R_{2}(E)<1.2551 \frac{E}{\ln E}-\frac{\frac{E}{2}-1}{\ln \left(\frac{E}{2}-1\right)} \\
& \frac{E}{4}-1.2551 \frac{E}{\ln E}+\frac{\frac{E}{2}-1}{\ln \left(\frac{E}{2}-1\right)}<R_{1}(E)<\frac{E}{4}-\frac{E}{\ln E}+1.2551 \frac{\frac{E}{2}}{\ln \frac{E}{2}}
\end{aligned}
\]

\begin{proof}
First of all, let $\frac{E}{2}$ be also an even number, then in view of the relation (8), by applying Dusart's inequality (1) on $\frac{E}{2}$ as
\[
\frac{\frac{E}{2}}{\ln \frac{E}{2}}-1<L_{2}\left(\frac{E}{2}\right)+R_{2}\left(\frac{E}{2}\right)<1.2251 \frac{\frac{E}{2}}{\ln \frac{E}{2}}-1
\]
and by relation (10) we thus arrive at the following inequality
\[
\frac{\frac{E}{2}}{\ln \frac{E}{2}}-1<L_{2}(E)<1.2551 \frac{\frac{E}{2}}{\ln \frac{E}{2}}-1
\]

But if $\frac{E}{2}$ is an odd number, then $\frac{E}{2}-1$ is an even number. Thus, in view of the relation (8), by applying Dusart's inequality (1) on $\frac{E}{2}-1$ as
\[
\frac{\frac{E}{2}-1}{\ln \left(\frac{E}{2}-1\right)}-1<L_{2}\left(\frac{E}{2}-1\right)+R_{2}\left(\frac{E}{2}-1\right)<1.2251 \frac{\frac{E}{2}-1}{\ln \left(\frac{E}{2}-1\right)}-1
\]
and by relation (11) we thus arrive at the following inequality
\[
\frac{\frac{E}{2}-1}{\ln \left(\frac{E}{2}-1\right)}-1<L_{2}(E)<1.2551 \frac{\frac{E}{2}-1}{\ln \left(\frac{E}{2}-1\right)}-1
\]

Since for all $x>141$ we have $\frac{\frac{2}{3}}{\ln \frac{2}{3}} \geq \frac{\frac{2}{3}-1}{\ln \left(\frac{2}{3}-1\right)}$, we thus for all even numbers $E>141$ obtain (16).

Since $\frac{E}{2}-L_{2}(E)=L_{1}(E)$, we can compute (17) by multiplying all sides of (16) with -1 and by adding $\frac{E}{4}$.

To get the relation (18), based on the relation (8), it is sufficient to subtract Dusart's inequality (1) applied on $E$ from the inequality (16).

To yield the inequality (19), based on the relation $\frac{E}{4}-R_{2}(E)=R_{1}(E)$, we follow the same reasoning that of (17).

At the end, it is necessary to mention that all inequalities are strict because the parameters $a_{E}, b_{E}, c_{E}, d_{E}$ always admit special values.	
\end{proof} 
In the equal of this section, By the previous Lemma we can deduce another inequality between the model elements, as follows:
\begin{lemma}\label{966}
For all even numbers $E>6$ we have $b_E<c_E$.
\end{lemma}
\begin{proof}
	Since for any real number $x>6$ we have $1.2551\frac{x}{\ln(x)-1}-\frac{x/2}{\ln(x/2)-1}<\frac{x/2}{\ln(x/2)-1}-1$, and by inequalities for $L_2(E)$ and $R_2(E)$ in the previous Lemma, the result is straightforward.  
\end{proof}
Ultimately, we can see the following order between the model elements, at least for any even number $E>22864$:
$$b_E<c_E<a_E.$$
\section*{3. Categorizing all Even Numbers into Distinct Structural Types by Using S.C.E Model}

In this section, we consider all possible structural types of representing an even number $E$ based on S.C.E model (2). To this end, since this model contains four elements, for a given even number $E$ we can then assume the following four general structures

\[
\begin{aligned}
	& d_{E} \leq \bigcirc \leq \bigcirc \leq \bigcirc \\
	& c_{E} \leq \bigcirc \leq \bigcirc \leq \bigcirc  \\
	& b_{E} \leq \bigcirc \leq \bigcirc \leq \bigcirc  \\
	& a_{E} \leq \bigcirc \leq \bigcirc \leq \bigcirc 
\end{aligned}
\]

where $\bigcirc$ denotes an arbitrary element of S.C.E model for $E$ other than the first element of each inequality. Therefore, we can expand each above general structure into 48 detailed structures which are listed in the following table
\newpage
\begin{center}
	\textbf{Table 1. Table of all typical structures of even numbers}
\end{center}
\begin{center}
	\includegraphics[width=1\textwidth]{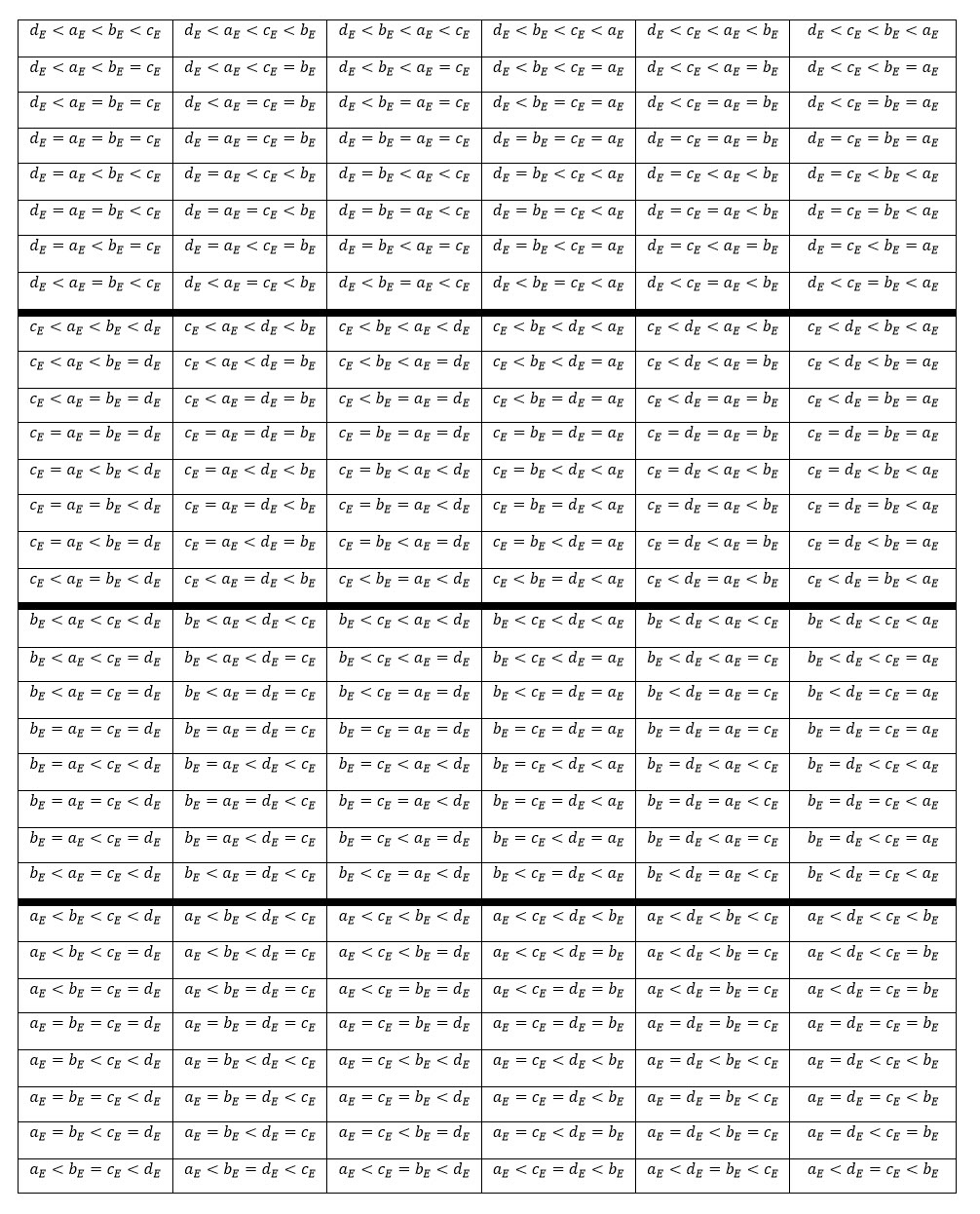}
\end{center}

\normalsize

If we remove all repeated structures to obtain distinct ones, we can see the following table which contains just 75 structures.

In fact after removing repeated structures, the first category to forth category contain $26,20,16$ and 13 structures respectively, which yield $75$ distinct typical structures based on S.C.E model. In the reminder of this paper, for all categories we prove that the element $d_{E}$ never vanishes. It is necessary to mention that the following three typical structures
	\newpage
\begin{center}
	\textbf{Table 2. Table of distinct typical structures of even numbers}
\end{center}
\begin{center}
	\includegraphics[width=1\textwidth]{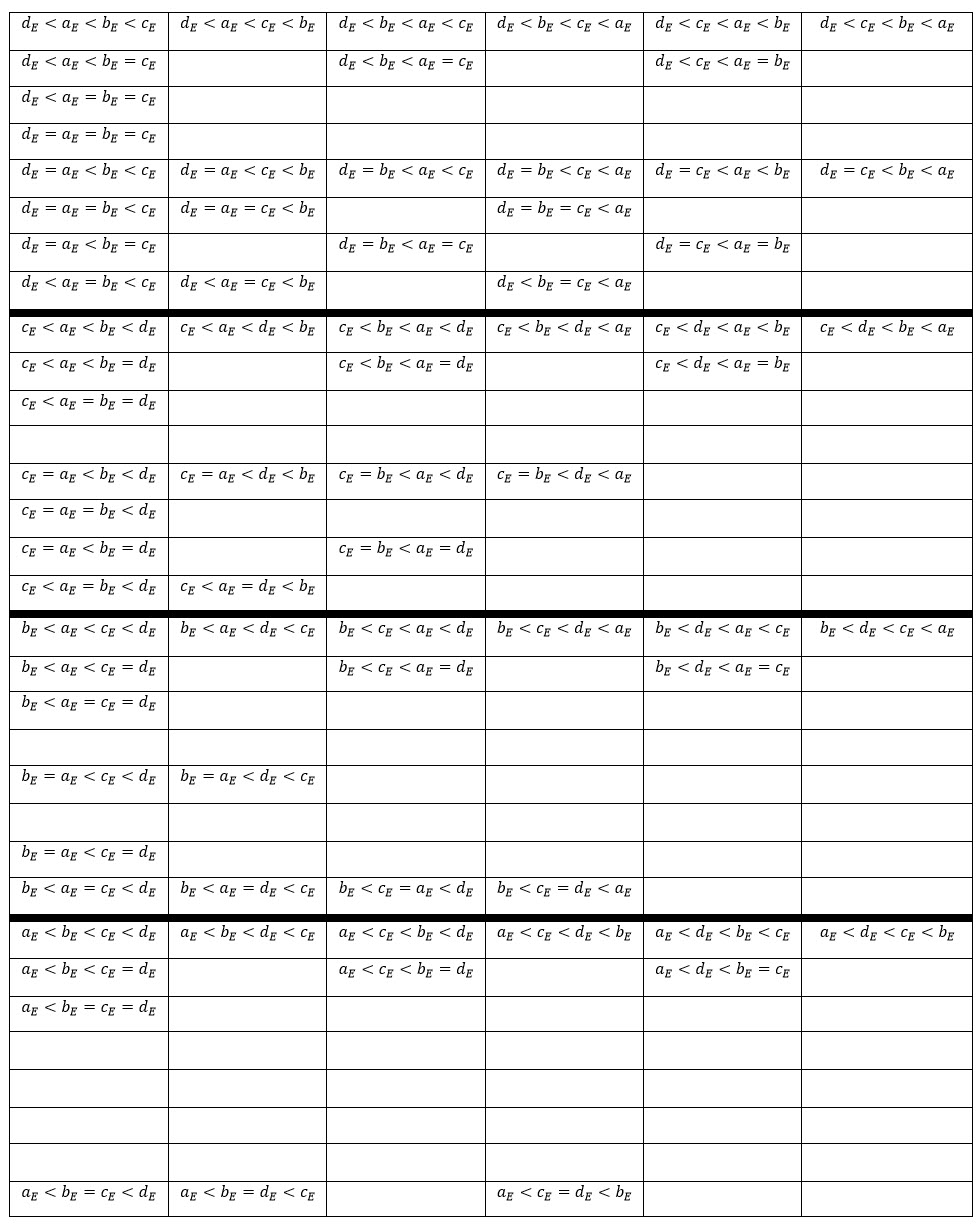}
\end{center}

cannot be proved to have nonzero element $d_{E}$

\[
\begin{aligned}
	& d_{E}<b_{E}<c_{E}<a_{E} \\
	& d_{E}<b_{E}=c_{E}<a_{E} \\
	& d_{E}<c_{E}<b_{E}<a_{E},
\end{aligned}
\]

where the last two inequality never happens by virtue of Lemma $\ref{966}$. For this reason, we prove Goldbach's strong conjecture relatively, namely, for 74 typical structures out of 75 ones. It is necessary to mention that, in this paper we come up with three inequalities using each of which we can knock out the last 75 typical structure. Unfortunately, the authors did their best to prove analytically these two inequality, but they did not succeed. On the other hand, we can see these inequalities hold true numerically.   

\section{Relative Proof of the Goldbach's strong Conjecture by Using S.C.E Model}

As we well know, the Goldbach's strong conjecture asserts that all positive even integers equal or greater than 4 can be expressed as the sum of two primes. It is necessary to mention that we prove this conjecture just for even numbers $E>2525$. To prove this conjecture, it is enough to show that for all positive even number $E>2525$ the element $d_{E}$ in the corresponding S.C.E model never vanishes. To this end, we first consider the following lemma.

\textbf{Lemma 4.1.} Let $E>130$ be a given even number and let $d_{E}=0$, then
\[
a_{E} \notin\left[0, \frac{E}{4}-1.2551 \frac{E}{\ln E}+1\right],\left[\frac{E}{4}-\frac{E}{\ln E}+1, \frac{E}{4}\right]
\]
\begin{proof}
This result is straightforward from the inequality $(\ref{2.1})$. 
\end{proof}
Now we can see the following main theorem.

\textbf{Theorem 4.2 (Goldbach).} Let $E>2525$ be a given even number and let corresponding $d_{E}$ be as in (2) and does not satisfy the inequality (24), then $d_{E}>0$.

\begin{proof}
First of all, since all elements of S.C.E model have nonnegative value, for the second category to the forth category in table (2) the element $d_{E}$ cannot be zero.

For the first category in table (2), we can take different policies. The first policy relates to the state that $d_{E}$ equals some other elements which consist of thirteen typical structures and can be listed as follows.
\begin{enumerate}
	\item 
	( $\left.d_{E}=a_{E}=0\right)$ : This item contains three typical structures and based on (15) it is not possible for $E>130$.\\
	\item 
	( $\left.d_{E}=b_{E}=0\right)$ : This item contains three typical structures and based on (18) it is not possible for $E>34$.\\
	\item 
	( $\left.d_{E}=c_{E}=0\right)$ : This item contains three typical structures and based on (19) it is not possible for $E>34$.\\
	\item 
	( $\left.d_{E}=a_{E}=b_{E}=0\right)$ : This item contains one typical structure and based on (16) and (19) it is not possible for $E>34$.\\
	\item 
	( $\left.d_{E}=a_{E}=c_{E}=0\right)$ : This item contains one typical structure and based on (17) and (18) it is not possible for $E>34$.\\
	\item 
	( $\left.d_{E}=b_{E}=c_{E}=0\right)$ : This item contains one typical structure and based on (18) and (19) it is not possible for $E>34$.\\
	\item 
	( $\left.d_{E}=a_{E}=b_{E}=c_{E}=0\right)$ : This item yields $E=0$ which is not the case.
	
	The second policy relates to applying the following two inequalities
	\[
	\begin{aligned}
		& a_{E}>c_{E}+2 d_{E}, \quad \forall E>2322.61 \\
		& a_{E}>b_{E}+2 d_{E}, \quad \forall E>2525.67
	\end{aligned}
	\]
	The inequality (25) is derived from the fact that for all $x>2322.61$ we have
	\[
	\frac{x}{4}-1.2251 \frac{x}{\ln x}+1>1.2551 \frac{x}{\ln x}-1
	\]
	and the inequalities (13) and (17). Similarly, the inequality (26) is derived from the fact that for all $x>2525.67$ we have
	\[
	\frac{\frac{x}{2}-1}{\ln \left(\frac{x}{2}-1\right)}-1>1.2551 \frac{x}{\ln x}-1
	\]
	and the inequalities (13) and (16). By using the inequalities (25) and (26), we see that for $E>2525.67$ the element $d_{E}$ cannot be zero in the following ten typical structures
\end{enumerate} 
	\begin{multicols}{2}
		\noindent
		$\displaystyle d_{E}<a_{E}<b_{E}<c_{E} \qquad d_{E}<a_{E}<b_{E}=c_{E}$\\
		$\displaystyle d_{E}<a_{E}=b_{E}=c_{E} \qquad d_{E}<a_{E}=b_{E}<c_{E}$\\
		$\displaystyle d_{E}<a_{E}<c_{E}<b_{E} \qquad d_{E}<a_{E}=c_{E}<b_{E}$\\
		$\displaystyle d_{E}<b_{E}<a_{E}<c_{E} \qquad d_{E}<b_{E}<a_{E}=c_{E}$\\
		$\displaystyle d_{E}<c_{E}<a_{E}<b_{E} \qquad d_{E}<c_{E}<a_{E}<b_{E}$
	\end{multicols}
	
	By counting individual items, we can see that 72 typical structures out of 75 ones are investigated in overall that their element $d_{E}$ never vanish and this completes the proof. 
\end{proof}
	Actually, we did our best to cover all 75 typical structures and in this way we can then prove Goldbach’s strong conjecture completely. We hope this can be done by using S.C.E model and this framework in the future. To this end, we invite all mathematician to take this framework into consideration to prove that remained typical structure of (24) cannot have zero value for $d_{E}$.
	
	In this way, we propose following three inequalities by using each of which we can pass by the last typical structure:
	\begin{itemize}
		\item 
		$b_{E}<\frac{E}{\ln E}-\frac{\frac{E}{2}}{\ln \frac{E}{2}}
		$, for each $10\leq E$
		\item 
		$c_{E}<\frac{\frac{E}{2}}{\ln \frac{E}{2}}-1$, for each $4\leq E$
		\item 
		$b_{E}+c_{E}<\frac{E}{\ln E}-1$, for each $2\leq E$
	\end{itemize} 
It is necessary to say that, we guess by following  
	\begin{center}
		
		\textbf{References}
	\end{center}
	
	\begin{enumerate}
		\item H. A. Helfgott. The ternary Goldbach conjecture is true. \textit{arXiv:1312.7748 [math.NT]}, 2013.
		\item P. Dusart. Explicit estimates of some functions over primes. \textit{Ramanujan J.}, 45(1):227--251, 2016.
	\end{enumerate}

\end{document}